\documentclass[12pt]{article}
\begin{document}
\bibliographystyle{plain}

\title
{Hot Hands, Streaks and Coin-flips: Numerical Nonsense in the New York Times}

%
%

\author{Dan Gusfield \\
Computer Science Department, University of California, Davis}



%
%

\maketitle

The existence of ``Hot Hands" and ``Streaks" in sports and gambling is hotly debated, but there is no 
uncertainty about the recent batting-average of the New York Times: it is now {\it two-for-two} in 
mangling and misunderstanding elementary concepts in
probability and statistics; and mixing up the key points in a recent paper that re-examines earlier
work on the statistics of streaks.
In so doing, it's high-visibility articles have added to the general-public's
confusion about probability, making it seem mysterious and paradoxical when it needn't be.
However, those articles make excellent case studies on how to get it wrong, and for discussions
in high-school and college classes focusing on quantitative reasoning, data analysis, probability and statistics. What I
have written here is intended for that audience.

\section{The Background}

The starting point for this discussion is an article by George Johnson in the New York Times Sunday Review on October
18, 2015, entitled ``Gambler, Scientists and Mysterious Hot Hand". That article discusses the claims in a recent
working paper (not yet peer reviewed) 
by two economists, Joshua Miller and  Adam Sanjurjo, entitled 
``Surprised by the Gambler's and Hot Hand Fallacies? A Truth in the Law of Small Numbers" \cite{MS2015}. According to the
Johnson article, the
Miller and Sanjurjo paper claims that the authors of a classic 1985 paper (Thomas Gilovich, Robert Vallone and Amos Tversky) \cite{GVT}
debunking the concept of hot hands in basketball, 
made an error in how they thought about probability. Quoting from the Johnson article:

\begin{quote}
A working paper published this summer has caused a stir by proposing that a classic body of research disproving the 
existence of the hot hand in basketball {\it is flawed by a subtle misperception about randomness}.  (italics added)
\end{quote}

Then, on October 27, 2015,
in a follow-on NYT article (in TheUpshot) entitled ``Streaks Like Daniel Murphey's Aren't Necessarily Random",
Binyamin Appelbaum wrote:
\begin{quote}
Last year two economists launched a more fundamental assault: They argued that disproofs of the ``hot hand" 
theory had {\it made a basic statistical error}. (italics added)
\end{quote}

Its a challenge to keep the players straight in this story, so to
recap, the issue of hot hands and the probability of streaks was first discussed in two 
academic papers and then in two subsequent NYT 
articles: the first paper, by Gilovich et al. in 1985 (which we will refer to as ``GVT"), claims
that the belief in hot-hands (in basketball) is not statistically supported; 
the second by Miller and Sanjurjo (which we will refer to as ``MS") in the summer of 2015 re-examines that work, suggests that 
statistical errors were made, 
and comes to a different conclusion;
the third by George Johnson in October 2015 in the NYT discussing the claims in MS; and a fourth article, by
Appelbaum ten days later in the NYT that repeats, and even strengthens, some of the statement from the Johnson article.

I am not interested in questions of hot-hands, streaks and gambling per se. 
Instead, my interest, and focus here, is
how the New York Times articles discus probability and statistics, and the confused and incorrect statements 
made in those articles. However, in order to explain the NYT errors, we will have to discuss streaks, hot hands,
and the two academic papers to some extent.

\section{The Central Technical Issues}

We want to identify the claimed ``subtle misperception about randomness" and ``basic statistical error" in GVT
that the two articles in the NYT are talking about. To do that, we have to say a bit about the statistical approach 
to the study of streaks and hot hands.

When trying to determine if streaks (successive baskets made,  heads on coin flips, wins in gambling, for example) have non-random
causes, such as skill or ``being in the zone",  the statistical
approach is to compare numerical features in observed data to features in data generated
at random. For example, suppose that a player makes a basket (a hit, coded as `H') on 50\% of the shots he takes, 
and that we have the entire record of the player's hits and misses.
We could look at that data
and ask what {\it percentage} of the Hs are
followed by another H. It has little effect in long sequences, but in a short sequence, we will compute the percentage by counting the
number of Hs in all but the last position, and the number of those Hs that are followed by another H (possibly in the last
position). For clarity, we give that percentage the name {\it HH-percentage}, although that term was not
used in the NYT articles, or in the academic papers. 
We could also determine the HH-percentage from data on Hs followed by a miss, coded as a `T'.
See Table \ref{HTtable}.

The HH-percentage might not be the ideal way to study questions of streaks and hot hands, although
a player with a few long streaks (who probably would be considered to have a hot hand) 
has a larger HH-percentage than a player with more, but shorter, streaks. Still,
the HH-percentage (in different terminology) is one of the first statistics examined in the GVT paper,
where they computed the HH-percentage for several individual NBA players in an individual season.
And, the HH-percentage is the {\it only} statistic that is discussed in the NYT articles, so it is the focus of
this note.
But, how specifically would we use the HH-percentage to determine if the player's Hs are unusually ``streaky", i.e., 
more concentrated into streaks that what we would expect by chance alone?  GVT says ``The player's performance, then,
can be compared to a sequence of hits and misses generated by tossing a coin."

Specifically, we could generate a long random sequence,
where each character in the sequence is
independently  chosen to be an H or a T with equal probability; and then compute the HH-percentage from that long sequence. 
We call that HH-percentage a ``reference number", and remember that it is obtained from a sequence that does not have 
any non-random influence.
Then, we would compare the HH-percentage obtained from the record of a chosen player to the reference number.
Intuitively, when a player's actual HH-percentage is computed from a long sequence (i.e., a large amount of data) it seems
appropriate to compare it to this reference number.\footnote{However, if
a player's HH-percentage is computed only from a ``relatively short" sequence (say a single game or even 
a single season), then the reference number
defined above might not be the most informative one to use.  Foreshadowing what will come later in this paper, this will be a key issue.}
If the reference number 
is very close to, or larger than, the HH-percentage in the player's record, 
then the player's HH-percentage does {\it not} support the conclusion that the player's streaks are due to
some non-random influence. That means, from the perspective of the HH-percentage, the player's baskets do 
not appear to be more streaky than do the Hs in 
a random sequence. Conversely, if the player's HH-percentage is ``significantly" larger than the reference value,
we {\it do} feel justified in thinking that some non-random influence is at work.
How much larger a player's HH-percentage must be in order to be ``significant", to 
support the assertion of non-randomness, is exactly the kind of issue that is
studied in statistics and probability theory, and is not our main concern here. 

Comparing a player's record to a randomly generated sequence is the basic statistical approach,
but do we actually need to generate a
random sequence in order to determine the reference value?  No. We might need
to generate random sequences to determine more complex statistics in random sequences, but
in the case of the HH-percentage, we don't need to generate any sequences 
because we {\it know} that the {\it probability}
of an H following an H is exactly the probability of an H on any individual flip, i.e, one-half.  
So, the observed HH-percentage  in a {\it long}
randomly-generated sequence will be about 50\%; about equal to the
frequency that an H is followed by a T, or a T is followed by an H.
That point should not be controversial or confusing.  

\subsubsection*{But the NYT article did confuse it}

Contrary to the point above, Johnson in the October 17 NYT article states:

\begin{quote}
For a 50 percent shooter, for example, the odds of making a basket are supposed to be no better after a hit -- still 50-50. 
{\it But in a purely random situation, according to the new analysis, a hit would be expected to be followed by another 
hit less than half the time.} (italics added)
\end{quote}

To be clear, the NYT article is talking about a ``purely random situation" of (memoryless) shots by a 50\% shooter, 
or equivalently, a sequence of fair coin flips. It is
not talking about some basketball-related phenomena (for example, a player being more tired or more closely guarded after making several
shots).  And, for even greater clarity, I interpret the statement `` ... in a purely random situation ... a hit would be expected
to be followed by a another hit less than half the time" as the same as ``... in a purely random situation ... the odds of making
a basket after a hit are less than 50-50. Equivalently, in a purely random situation ... the probability that
a hit will be followed by another hit is less than one-half."\footnote{If you think this interpretation is wrong, then you will probably
find the rest of this paper wrong, and can stop reading now.}

\section{Really!?!}
Can that statement about hits (and coin flips) in Johnson's article be correct, that ``in a purely random situation ... 
a hit is expected to be followed
by another hit less than half the time?" 
Surely, there is something wrong here, because in a purely random situation
{\it every} flip will be an H with the same probability that it is a T --- exactly one-half. So, a hit is {\it expected} 
to be followed by another hit (H) one-half of the time, which is as often as it is expected to be followed by a miss (T).
Several of the
on-line comments to the NYT submitted by readers after the publication of Johnson's article 
correctly pointed this out, and even identified the source of Johnson's confusion, which we will discuss in detail below.
But, despite the readers comments, ten days 
later, Appelbaum in the NYT article (in TheUpshot), doubled down on Johnson's statements, making even more explicit
statements:

\begin{quote}
Flip a coin, and there's an equal chance it will land heads or tails. Researchers had treated that 
50 percent chance as the definition of a random outcome. But Joshua Miller of Bocconi University and Adam Sanjurjo of the 
Universidad de Alicante pointed out something surprising: {\it In the average series of four coin flips, the sequence heads-heads is 
significantly less common than heads-tails.} (italics added)
\end{quote}

Really? In the table of coin flips (similar to Table 1 below) that Appelbaum directs the reader to examine, 
heads-heads occurs {\it exactly} the same number of times that heads-tails occurs. So is Appelbaum's statement
pure nonsense, or is it based on some truth, but one that is very poorly stated?  He continues:

\begin{quote}
On average, just 40.5 percent of the heads are followed by another heads. 
Yes, this sounds crazy. But it happens to be true.\footnote{One might argue that Appelbaum has 
a tiny, tiny bit a wiggle room, because he
does not define what ``the average series of four coin flips" means, or what ``on average" means in the second quote.  
But he directs the reader to the Johnson article with
the table 
showing that exactly 50\% of the heads, in the first three positions, are followed by another head. 
So, his statement is particularly confused and incorrect.}
\end{quote}

And, this assertion has consequences for the study of streaks. Referring to MS, Appelbaum writes:

\begin{quote}
The implication, they argued, is that past studies had set the bar too high. Streaks that has looked like random luck
were actually statistically unlikely. The ``hot hands fallacy", they wrote, was remarkably persistent because it was true.
\end{quote}

\section{So What is Going on?}

Both NYT articles imply that the ``basic statistical error" made in GVT is to assume that
the probability that an H will follow an H is one-half, 
in sequences of heads and tails created by flips of a fair coin. In the case of the sixteen length-four sequences, the
implied ``error" in GVT is
the assumption that 50\% of the heads that are followed by another flip, are followed by  another heads.
But they are. So what is going on here?

\medskip

{\bf Spoiler alert}: In trying to interpret MS, the Johnson article made incorrect and imprecise statements about probability and statistics.
The Appelbaum article repeated, more strongly, the main one. Both papers miss the key points made in MS.
In truth, in purely fair coin flips, each H (other than the last one in the sequence) 
will be followed by another H {\it with probability one-half}. Period. 
Miller and Sanjurjo also make that clear. So how did Johnson and Appelbaum get it so wrong? 

\medskip

\subsection{The Johnson Table}
Following a similar example and table 
in the MS paper (but not a similar conclusion), here is what Johnson did in his article. 
He looked at the sixteen, length-four sequences shown in Table \ref{HTtable}. For each sequence that contains an H in
one of the first three
positions (there are fourteen of these) he calculated the {\it percentage} of those Hs that are followed by another H.
For example, in the sequence HHTT, the percentage is 50\%, and in HHHH it is 100\%, and in HTTT it is 0\%. Hence,
Johnson calculated the individual HH-percentage for {\it each} of the relevant fourteen sequences. Then he added those fourteen
HH-percentages,
divided by fourteen, and got about 40.5\%. That is, he {\it averaged} the HH-percentages calculated from the fourteen 
relevant sequences. As he writes:

\begin{quote}
... calculate for each sequence the odds that a head is followed by a head and average the results. The answer is not
50-50, as most people would expect, but 40.5 percent -- in favor of tails.
\end{quote}

All true. The arithmetic is right, and the 40.5\% average may indeed seem surprising to some people. But so what?
What does that average have to do with the probability that an H 
is followed by another H? 
Nothing! It is nonsense to conclude from that averaging 
that ``a hit is expected to be followed by a  hit less than
50\% of the time", or that ``On average, just 40.5 percent of the heads are followed by another heads."

\medskip

\subsection{Counting Bathrooms}
In order to explain what Johnson and Appelbaum got wrong, we 
look here at a more extreme scenario. Suppose we want to calculate the average number of bathrooms in 
the houses in the U.S.  The right way to calculate this is to find the number of bathrooms in each of the (millions) of U.S. houses,
sum up those numbers and divide by the number of houses in the U.S. But here is another suggestion: After finding the
number of bathrooms in each of the houses, divide the houses into two groups: those that have more than 30 bathrooms, and those that
have 30 or fewer bathrooms. (San Simeon, the former country house of the Hearst family, has 61 bathrooms, and the White House has
35). Next, compute the average 
number of bathrooms in the 
first group of houses (perhaps that average is 32.5 bathrooms), and compute the average for the second group of houses 
(around 2.7 in a recent survey).
Finally, average those two averages, to get 17.5 bathrooms (just slightly more than I have in my house).  And even though
I made up the average of 32.5 for the first group, the correct average in the first group will be at least 30 (why?), 
and the average in the second group
is actually close to 2.7, so the true average of those averages will be larger than 16.
Probably (but what does that really mean?), the average of 16 bathrooms per U.S. house does not mesh with your sense of reality. 
So what went wrong?

By averaging the two averages, we give {\it equal} weight to each of the averages, ignoring the fact that the first average
comes from a very small number of houses, while the second average comes from a huge number of houses. That kind of
average is called an {\it unweighted} average. But, to get the
correct average number of bathrooms, you must give equal weight to each {\it house}, not to each {\it group of houses}.

Now if for some reason you don't have data on the number of bathrooms in each individual house, but are given the two averages 
in the two groups, and are also given the number of houses in the two groups, you could
multiply the first average by the number of houses in that group,
multiply the second average by the number of houses in that group, and add the two products to get the total number
of bathrooms in the U.S. Then, to get the correct average number of bathrooms, you would divide that total by the sum of
the number of houses in the two groups, i.e., the total 
number of houses. This is called a {\it weighted average} of the averages, and would give a result of about 2.7. 
Note that computing the weighted average  is just a backwards way of doing what we would do to compute the average number of
bathrooms in a U.S. house, if we
had the raw data on each house: find the total number of bathrooms and divide by the number of houses.

\paragraph*{Back to HH-percentages}
How does the bathroom story relate to HH-percentages? There are 24 Hs that occur in the first
three positions of the 16 sequences of length four. These 24 Hs are analogous to the houses in the bathroom story.
If you want to
compute the percentage of those 24 Hs that are followed by an H, or equivalently, how often ``a hit is expected to be
followed by a hit",
you should {\it not} divide those 24 Hs into groups (in this case, 14 groups, each called a ``sequence"), 
find the HH-percentage in each group (sequence), and then average those percentages. To do so gives equal weight to each group 
(sequence), ignoring the fact that some groups (sequences) have more Hs than others do. That is, you should
not compute an {\it unweighted} average of
the HH-percentages. Instead, to calculate the probability that an H follows an H, you need to give equal weight to
each H that occurs in the first three positions of some sequence, or if you start from the HH-percentages of the fourteen sequences, you need
to compute a {\it weighted} average of those HH-percentages; each HH-percentage 
weighted by (multiplied by) the number of Hs in the first three
positions of the sequence that the HH-percentage comes from.

Another numerical reflection of the difference between unweighted and weighted average HH-percentages is the fact that
in the sixteen length-four sequences, there are only {\it eight} that have any occurrence of HH, but there are {\it eleven} that have an 
occurrence of HT. That is, the distribution of HH and HT is not uniform in the fourteen sequences.\footnote{At first, this may seem paradoxical since the two counts might be expected to be equal by ``symmetry". But,
the two occurrences are {\it not} symmetric, which I leave you to ponder.}   Similarly, there are seven sequences that have HHH, but eight that
have HHT. So, in random sequences, {\it if} your unit of analysis is the whole sequence, you will observe a T following an H more often
(in more sequences) than an H following an H.  You will also observe an H following a T in more sequences than an H following an H.
So, by equally weighting the sequences, we under-represent the HHs and over-represent the HTs.

\paragraph*{The Take-Home Lesson:} The 
unweighted average of the averages calculated from non-overlapping subsets of a set is 
not always equal 
to the average in the entire set.
That is just a numerical fact, and is elementary text-book material in
any basic statistics book or course.
The numerical example in Johnson's article does nothing more than illustrate that 
fact in the case of all possible length-four sequences of 
fair coin flips.  It does {\it not} establish that
``In a purely random situation ... a hit would be expected to be followed by another hit less than half the time."

\subsection{The MS Table and Unweighted Averages} 
While Johnson and Appelbaum completely miss the issue of weighted versus unweighted averaging, 
Miller and Sanjurjo understand it perfectly well, as did many of the NYT readers who commented on the Johnson and Appelbaum articles. 
MS contains a table that is similar to the one in Johnson's article (and to Table 1 below), and it obtains the same average, but MS
does not state the conclusion
that Johnson and Appelbaum do. 
In fact, in a blog discussion this summer, Miller states:

\begin{quote}
We do not assert that: ``a way to determine the probability of a heads following a heads in a fixed sequence, you may calculate the 
proportion of times a head is followed by a head for each possible sequence and then compute the average 
proportion, giving each sequence an equal weighting"  ... 

it is a mistaken intuition to treat this computation as an unbiased estimator of the true probability.
\end{quote}

MS begins by stating that if one million fair coins are each flipped four times, and an 
HH-percentage\footnote{The actual terms they use are ``relative frequency" and ``empirical probability".} 
is obtained for each coin, those million
HH-percentages would average to ``approximately 0.4".  In explaining this, they state:

\begin{quote}
The key ... is that it is not the flip that is treated as the {\it unit of analysis}, but rather the {\it sequence} of flips from each
coin ... (italics added)

Therefore, in treating the sequence as the unit of analysis, the average empirical probability across coins amounts
to an unweighted average\footnote{For clarity, note that in this quote, it is implied that ``sequence of flips"  is ``sequence of four flips".} ...
\end{quote}

The unweighted average of averages (about 0.405) is not equal to the probability (exactly 0.5) 
of an H following an H in
four fair coin flips.  
The NYT articles printed nonsense, because what they wrote 
suggests that 
these are the same.\footnote{The Johnson paper is actually more confused and confusing, because, as explained above, 
it suggests that the probability that an H follows an H
is less than half, and yet it also points out that in the sixteen sequences of length four, the number of Hs that are followed by another
H is exactly the same as the number of Hs that are followed by a T. It tries to explain this apparent contradiction by introducing the
concept of a ``selection bias". This is actually more nonsense; we will return to this later.}

\paragraph*{But why?}
The table in the Johnson article, which the NYT articles misunderstand,  originates in the MS paper. But why? 
One of the reasons MS examines {\it un}weighted averages is explained next.

\section{Modeling the Gambler's Fallacy}

The Miller and Surjurjo  paper is concerned with several streak phenomena {\it in addition} to hot hands in sports.
The main one is called the ``Gambler's Fallacy", which is the belief 
that a streak (winning or losing) in a game of pure (or mostly pure) chance,
will soon be reversed, in order to achieve the {\it long-run} expected win/loss frequency.\footnote{In the ``long run", the frequency of wins should
be about equal to the frequency of losses. That is a consequence of the ``law of large numbers". The belief that we should also see this balance in
small sequences has been facetiously called the ``law of small numbers".}

This fallacy is most clearly defined in terms of a sequence of fair coin flips, where the gambler's fallacy is: 
\begin{quote}
If one observes a growing streak of Hs, the probability
that the next flip will be T increases after each successive H. That is, the longer the streak of Hs, the higher
is the probability that the next flip will be an T.
\end{quote}

Restricted to just two consecutive flips, the gambler's fallacy is that the probability that an H will be followed by another H
is lower than the probability that it will be followed by a T.
Thus, the gambler's fallacy is similar to the belief in a ``hot hand", but there an H is believed to be {\it more} likely, rather than
less likely, after an H.

Both GVT and MS assert that the gambler's fallacy is a commonly
held belief. Of course, this belief is a fallacy, since 
the probability that the next flip will be an H is precisely one-half (in a fair coin),
no matter what the past history is.  

MS uses {\it unweighted} averages of HH-percentages, because 
Miller and Sanjurjo assert that peoples' {\it beliefs} about streaks  in gambling
are based on gamblers' observations of many {\it short, but whole} sequences of events, or {\it complete games}. 
These are the ``units of analysis" 
that best {\it model} how people incorrectly come to believe in the gambler's fallacy.  
In the analogy of coin flips, 
a finite sequence of flips (say, of length four) is the unit of analysis, and multiple sequences are observed. 
Miller and Sanjurjo assert
that people use ``natural" statistics, which
{\it equally} weight what they observe in each sequence or game. Hence, their beliefs
are essentially based on an {\it un}weighted averaging of the sequences and games they observe. And since
unweighted averages of HH-percentages underestimate the true probability that an H will follow an H,\footnote{The unweighted
average is 40.5\% for sequences of length four. For longer 
sequences, the unweighted average HH-percentage remains less than 50\%, although it approaches 50\% as the
sequence length increases. For example, 
in length-six sequences, the unweighted average HH-percentage is
41.6\%, averaged over the 62 sequences that have an H in one of the first five positions.}
this intuitive (but incorrect) 
thinking leads to a belief in the gambler's fallacy. Miller and Sanjurjo write:

\begin{quote}
The implications for learning are stark: to the extent that decision makers update their beliefs regarding
sequential dependence\footnote{The term ``sequential dependence"  refers to the way that one event relates to
a prior one. In the case of two flips, it refers to whether an H or a T follows an H.}
with the (unweighted) empirical probabilities that they observe in finite length sequences, they can never unlearn
a belief in the gambler's fallacy. ...

no amount of exposure to these sequences can make a belief in the gambler’s fallacy go away. 
\end{quote}

And:
\begin{quote}
... in treating the sequence as the unit of analysis, the average
empirical probability across coins amounts to an unweighted average ...
and thus leads the data to appear consistent with the gambler's fallacy.\footnote{The phrase `across coins' should be
interpreted as `across sequences' in the treatment here, because in MS it is assumed that each sequence is generated
by a fair, but different coin.}
\end{quote}

\section{But What About Basketball?}

We have seen why MS is concerned with unweighted averages of HH-percentages in their treatment of the Gambler's Fallacy.
But what about streaks in basketball? MS is concerned with {\it un}weighted averages there also, but the
explanation for this is more subtle than for the Gambler's Fallacy. To get to that explanation, we first have to
discuss another way that Miller and Sanjurjo explain their main statistical observation.

\subsection{Alice and Bob}

In trying to explain the main technical issue in their paper, Miller (in an online post) describes a competition
between two players I will call Alice and Bob (I am modifying the description of the game, but not altering its
mathematical features).

The scenario is as follows: A computer has been programmed\footnote{Both Alice and Bob have previously verified that 
the program is correct.} to simulate a fair coin flip.
It first generates
a random sequence of four fair coin flips (printing out the sequence for later verification, and Bob can't see the output now); 
then, the computer randomly picks a position of one of the Hs in the sequence, 
provided that it is in one of the first three positions.
If there is no such position, the computer starts again.
If there is such a position, Bob is invited to bet whether the following position in the sequence is an H or a T. Note that
the value (H or T) has already been generated and written down. If Bob's bet is correct, 
Alice pays
him \$1, and if it is incorrect, Bob pays Alice \$1.
Notice that Alice has no active role except to pay out or collect the winnings. What should Bob pick, H or T?
One is tempted to answer ``pick either one, because on any flip the probability of an H is the same as the probability of a T."
But that answer ignores the full context of the competition.

The answer is that Bob should pick T, not H.  Randomly generating a sequence of length four is equivalent to randomly picking one of
the sixteen sequences shown in Table 1, because the probability of generating any specific sequence is that same as the probability
of generating any other sequence (i.e., $(\frac{1}{2})^4$). So, instead of imagining the computer generating a random sequence of length four, imagine that
the computer randomly (with equal probability) picks one of the fourteen relevant sequences; and then randomly picks an H in one of the
first three positions of that sequence, at which point Bob bets either that the following (already determined) flip is an H or is a T.
If we repeated this scenario many times, the frequency
that the following character is an H, would be a good estimate of the {\it unweighted} average of the HH-percentages, 
over the fourteen relevant sequences. 
Since the unweighted average of the HH-percentages is 40.5\%,
the probability that Bob will win if he picks H is only 0.405.  That is why Bob should pick T.

The key point is that this scenario has {\it two} stages: the computer {\it first} picks a sequence with equal probability; and {\it second}
it randomly picks an
H in the first three positions of that sequence (if there is one). But that is very different from a {\it one-stage} scenario
where the computer randomly picks an 
H in one of the first three positions of the fourteen relevant sequences. In this second scenario, the sequences 
would not be picked with equal
probability, because the distribution of Hs is not uniform. In the second scenario, the probability that Bob wins if he picks H 
{\it is} exactly 0.5, not 0.405.
The first scenario corresponds to an unweighted averaging of the HH-percentage observed in each sequence, 
and the second corresponds to a weighted average of the HH-percentages. Further, the second scenario roughly reflects how GVT
obtained the reference number it used to compare a player's HH-percentage, while the first one roughly reflects how MS 
does.\footnote{In this paper, we have only discussed HH-percentages because it is the only statistic discussed in the NYT articles,
and it is sufficient to illustrate the key difference in the approaches of GVT and MS.  But actually, the main statistic
discussed in both GVT and MS is a bit more involved. 
Define the ``TH-percentage" in a sequence as the percentage of Ts that are followed by an H.
Then define statistic $D$ for a sequence as its HH-percentage {\it minus} its TH-percentage.
$D$ relates the relative frequency that an H follows an H to the relative frequency that it follows a T in a sequence, and it
may be a more meaningful statistic to use to answer questions about ``hot hands".
Now, considering again all the sequences of length four, we see that there are exactly the same number of HH pairs as TH pairs.
However, the unweighted average of the sixteen $D$ values is not zero, but something {\it less} than zero. This is analogous
to the fact that the unweighted average HH-percentage is less than 50\%, the percentage of all Hs in the first three positions that
are followed by another H. So, the issues that arise in using $D$ values are well illustrated
by considering only the HH-percentages.} 

\subsection{Is this relevant for basketball?}

The answer depends on what specific question you are asking.
For example, we could ask: 
Did a specific player exhibit
``streak shooting" in a specific game? or ask:
Is a specific player a ``streak shooter" generally, considered over a season or
their entire career?

\subsubsection{Analysis for a single game}
For the first question,
let's suppose that a player, who has a long-term 50\% hit rate, shoots four times in a game. We want to know if
the player exhibited a hot hand, and so we compute his HH-percentage for that game. We  then
compare that number to a reference number derived from random sequences generated without an explicit hot hand.
We could compare his number to the HH-percentage generated from a {\it long} random sequence, in which case the reference
number should be 
50\%. This essentially (but not exactly) reflects the approach in GVT.
But another approach, which reflects a different way to {\it model} a player without a hot hand, is to consider
all the random sequences of length four.
If we use that model, then the number to compare with is 40.5\%. The reasoning, detailed next, is similar to the
reason that Bob should bet on T rather than H. 

We model the player {\it without a hot hand} simply
as a fair coin, i.e., each shot is a hit (H) with probability of one-half, independent of any other shot;  
so in a game, the player (with no hot hand)  generates
a random sequence of Hs and Ts. 
As discussed earlier,
we can also think of the generation of a random sequence of four flips as a random {\it selection} (with equal probability) 
of one of the sixteen four-flip sequences. 
Thus, instead of thinking of the player (who does not have a hot hand) {\it generating} a random sequence of length four, we
model the player's record as a selection of one of the sixteen sequences, chosen at random.\footnote{But, since we 
are only interested in the hot-hand question, the only relevant sequences considered in MS are the 
fourteen sequences that have an H in one of the first three positions. I would have chosen to include
the other two sequences as well, on the grounds that a player who makes none of his shots, or only his last shot,
should certainly not be said to have a "hot hand".} 
So, we compare the player's actual HH-percentage in the game
to the HH-percentages from the fourteen relevant random sequences. But which specific number obtained from those
HH-percentages should we use? The statistical approach is to consider what we would see over time, if we randomly selected many sequences of
length four from the relevant fourteen.
Each random sequence is selected with the same 
probability, 1/14, and so if we select many random sequences of length four and take the average of the HH-percentages we observe, 
what we will get is the sum of the fourteen HH-percentages, divided by 14, i.e., the {\it unweighted average} of the fourteen HH-percentages. 
So, in this model of a player without a hot hand, we should compare the real player's HH-percentage to 40.5\%.

This means that when the unit of analysis is an individual sequence, rather than an individual flip, to determine if
a player with a 50\% hit rate exhibits streak shooting in a single sequence (a game, say), 
we should not compare the observed  HH-percentage
to 50\%, but rather to a number less than 50\%.

\paragraph*{Four is just for illustration}
Now, in most games, a player shoots more than four times, and in fact, shoots a different number of times in 
each game. So the four-shot example is just an illustration, a simple idealized scenario used to explain the point made in MS:
when the unit of analysis is a player's 
record in a {\it single} game, or perhaps even a single season, 
the value we compare to should be lower than the long-term hit rate of that player. 
So, for example, if we observe that a player with a well-established hit-rate of 50\% has an HH-percentage of 50\%
in a game or season, that 
can be taken as evidence that the player {\it has} exhibited a hot hand in that sequence, rather than evidence against it. 
How strong that evidence is in favor of a hot hand requires additional probability theory, and is affected by the length of the
sequence. Length four sequences demonstrate the effect dramatically, but over a season or career, the sequence might be long enough that
the effect is small.
For a fair coin, the average HH-percentage, averaged over all sequences of length $k$, approaches one-half as $k$ increases, 
although it is always below one-half.

\subsubsection{Season or career-long analysis}
For the second question, if a career is long enough and an individual player makes many shots, it
seems appropriate to compute an HH-percentage over all of the Hs, equally weighting each basket, meaning that
the unit of analysis is an individual H. This might even be sensible for a single season, depending on the number of shots taken.
That is essentially the approach taken in GVT, where a player's HH-percentage  is compared to his season-long hit rate.
If they are close to each other, then GVT takes that as evidence against a hot hand.  

But according to MS, in a random sequence of coin flips with length equal to the number of shots, call it $K$, that a player makes in
a typical season, the unweighted HH-percentage averaged over all the possible $K$-length sequences,
is still significantly less than 50\%. Hence, MS assert that the proper unit of analysis is a player's record for an 
entire season, considered as
one sequence. In that case, when determining if a player (with a 50\% hit rate) was generally a streak shooter, we should
compare his HH-percentage for the season to  
the unweighted average HH-percentage in all the possible $K$-length sequences. Then,
as in our discussion of a single game, 
an HH-percentage of 50\% for the season should be taken as
evidence that the player is a streak shooter.

\subsection{It's the model}
When there is a dispute between academics, particularly in science or mathematics, it is easier for a journalist to 
explain the dispute by saying that one of the parties made an ``error" or had a ``misperception". And, that explanation
may be more attractive to the public. But the reality is
often that the parties have a legitimate difference of opinion on some methodological or data issue.
When using mathematics to study a natural or human phenomenon, we must create a detailed
{\it model} of the phenomenon to allow the application of mathematics. Different models can lead to different
ways that mathematics is used. In the disagreement between the GVT and MS papers, the fundamental issue is not that one of the parties made
a mathematical error or had a misperception of randomness --- the underlying issue concerns the ``unit of analysis" that the
mathematics applies to, and that is determined by the way one {\it models} a player without a hot hand.
The unit of analysis then dictates whether an unweighted or weighted average of HH-percentages (over all random sequences of a fixed length) 
is used to determine the reference number that a player's HH-percentage will be compared to. 
The take-home lesson here is that modeling is a critical and difficult part of the application of mathematics.
It is not enough to ``get the math right". To make the math meaningful, you have to create a meaningful model, and
people often disagree on which models are the most meaningful.

\section{One Last Piece of Nonsense}
In the caption of the table shown in Johnson's article, after showing that the computed average is 40.5 percent, Johnson
adds:
\begin{quote}
This is not, however, a violation of the laws of randomness. A head is followed by a head 12 times and by a tail 12 times.
{\it But by concentrating only on the flips that follow heads and ignoring the other data, we are fooled by a selection bias.}
(italics added)
\end{quote}

What? The disagreement between the 40.5 percent average, and the fact that in the table a head follows a head
exactly the same number of times that a tail follows a head, has nothing to do with ``concentrating only on the flips
that follow heads." A ``selection bias" is discussed in MS, but it is
the consequence of choosing one of the fourteen relevant sequences of length four, with equal probability, independent of 
how many Hs are in the first three positions.
As we have discussed above, since there  are more sequences that have at least one HT than have at least one HH,
the selection bias leads to seeing a T follow an H more often than an H following an H. It is nonsense
to say that ``by concentrating only on the flips that follow heads and ignoring the other data, we are fooled by
a selection bias."

\section{Aristotle and Appelbaum}

Appelbaum, after asserting that ``On average, just 40.5 percent of the heads are followed by another heads" continues with

\begin{quote}
Go ahead, see for yourself (link to the October 17 NYT article by George Johnson).
\end{quote}

That link leads to the Johnson NYT article, which contains the table 
showing that precisely 50 percent of the heads
are followed by another heads. So, although Appelbaum encourages the reader to ``see for yourself", it seems that he did not
make the effort. Apparently, he was so convinced of the claim that he didn't think it needed empirical testing.

This reminds me of the story about Aristotle and the role of theoretical versus empirical thinking.
Aristotle asserted that men have more teeth than women. As Bertrand Russell wrote: ``Aristotle could have avoided the mistake
of thinking that women have fewer teeth than men, by the simple device of asking Mrs. Aristotle to keep her mouth open while
he counted." 

So Appelbaum didn't do the counting. 
Ironically, both the NYT articles discuss the the psychology of perceived randomness, and
how easy it is to be fooled, even in the face of clear evidence.
Johnson writes:

\begin{quote} 
For all their care to be objective, scientists are as prone as anyone to valuing data that supports their hypothesis over those
that contradict it.
\end{quote}

\section{What about the New York Times, the Newspaper of Record}

How could the Johnson article, and even more the Appelbaum article, have been published in the New York Times? They wrote
nonsensical things about probability and seriously misunderstand MS. The Wall Street Journal also wrote
about the hot hands dispute in ``The `hot hand' debate gets flipped on its head", by Ben Cohen, September 29, 2015, and
initially made exactly the same mistake as the NYT articles. They wrote:

\begin{quote}
Toss a coin four times. Write down what happened. Repeat that process one million times. What percentage of flips after heads also
come up heads? The obvious answer is 50\%. That answer is also wrong. The real answer is 40\%...
\end{quote}

But then on September 30, in an online version of the article, the error is noted and corrected to:
\begin{quote}
Toss a coin four times. Write down the percentage of heads on the flips coming immediately after heads. 
Repeat that process one million times. On average, what is that percentage?

...

{\bf Corrections \& Amplifications:}

A previous version of this article incorrectly describes the question regarding coin flips. The question is about the 
average percentage of flips, not the overall percentage of flips. (Sept. 30)
\end{quote}

But the NYT did not make any correction of the mistakes in the Johnson article.
More embarrassingly, since several of the
readers of the Johnson article correctly pointed out the nonsense, how did the Appelbaum article make it past the editors? 
As one of the readers (Larry from St. Louis) commented online after the Appelbaum article:

\begin{quote}
It is shocking that such a basic error would get through both the Sunday Review and the Upshot. Is no one on the 
paper paying attention to what people write? Further, if the author of this Upshot column read the comments on the Sunday Review article, then he would have figured out the error for himself.
\end{quote}

Apparently
neither Appelbaum nor the editors read the comments of the readers, or if they did, they didn't understand them, 
or think to ask an expert.
And now, almost two months after the publication of the Johnson and Appelbaum articles, and in contrast to the WSJ, 
there is no retraction, or further clarification in the 
Times, or even the printing of a letter to the editor.  As an educator in a field involving mathematical reasoning,
and one concerned with the public's understanding of quantitative issues, and a long-time NYT subscriber\footnote{and not a WSJ reader}, 
this is all very disturbing.

\begin{table}
\begin{center}
\begin{tabular}{lcccr} 
The 16 & Number of Hs in the &  Number & HH-Percentage  & \\
sequences & first three positions & of HHs &   &\\
 & &  &   &\\
\hline 
HHHH & 3 & 3 & 100 & \\
HHHT & 3 & 2 & 66.66 ... & \\
HHTH & 2 & 1 & 50 & \\
HHTT & 2 & 1 & 50 & \\
HTHH & 2 & 1 & 50 & \\
HTHT & 2 & 0 &  0 & \\
HTTH & 1 & 0 &  0 & \\
HTTT & 1 & 0 &  0 & \\
THHH & 2 & 2 & 100 & \\
THHT & 2 & 1 &  50 & \\
THTH & 1 & 0 &  0 & \\
THTT & 1 & 0 &  0 & \\
TTHH & 1 & 1 &  100 & \\
TTHT & 1 & 0 & 0 & \\
TTTH & 0 & 0 & 0 & \\
TTTT & 0 & 0 & 0 &  \\
\hline
Total & 24 & 12 &  &\\
Average from &    &    & 40.5 &\\
the first 14 & & & & \\
sequences & & & & \\
\hline
\end{tabular}
\end{center}
\caption{The sixteen HT sequences of length four. The first fourteen contain an H in the first three positions. 
In each
of those fourteen sequences, the number of Hs in the first three positions is shown; next the number of those Hs that
are followed by another H is shown; then the percentage of Hs in the first three positions that are followed 
by another H is shown. This is the HH-percentage.  The total number of Hs in the first three positions is 24, and the 
number of Hs in the first three positions that are followed by another H is 12,  exactly 50\%. However, the unweighted average of
the percentages is not 50\%, it is about 40.5\%. True, but so what?}
\label{HTtable}
\end{table}

\end{document}